\documentclass[a4paper,11pt]{amsart}
\usepackage{amsfonts,amssymb,epsfig,latexsym}
\usepackage{amsmath,comment,bm,mathrsfs,braket}
\usepackage[latin1]{inputenc}
\usepackage{color,layout}
\usepackage{ifthen}

\usepackage{xparse}
\makeatletter
\RenewDocumentCommand{\title}{om}{%
   \IfNoValueTF{#1}
     {\gdef\shorttitle{Resonance free domain}}%
     {\gdef\shorttitle{#1}}%
   \gdef\@title{#2}%
}
\makeatother
\usepackage{graphicx}
\usepackage{color}  
\usepackage{bmpsize}

\newtheorem{theorem}{Theorem}[section]

\newtheorem{remark}[theorem]{Remark}

\def\square{\hbox{\vrule\vbox{\hrule\phantom{o}\hrule}\vrule}}

\newcommand{\be}{\begin{equation}}
\newcommand{\ee}{\end{equation}}

\newcommand{\til}[1]{\widetilde{#1}}

\numberwithin{equation}{section}

\newcommand{\N}{\mathbb{N}}

\newcommand{\R}{\mathbb{R}}
\newcommand{\C}{\mathbb{C}}

\newcommand{\cS}{{\mathcal S}}

\newcommand{\cB}{{\mathcal B}}

\newcommand{\p}{\partial}

\newcommand{\e}{\varepsilon}
\newcommand{\dl}{\delta}

\newcommand{\re}{{\rm Re}\hskip 1pt }
\newcommand{\im}{{\rm Im}\hskip 1pt }
\newcommand{\ord}{{\mathcal O}}

\newcommand{\wey}[1]{#1^w}
\newcommand{\ope}[1]{{\operatorname{#1}}}
\newcommand{\mc}[1]{{\mathcal{#1}}}

\numberwithin{equation}{section}

\begin{document}

\title{Resonance free domain for a system of Schr\"odinger operators with energy-level crossings}
\author{K.~Higuchi}


\maketitle

\begin{abstract}
We consider a $2\times 2$ system of 1D semiclassical differential operators with two Schr\"odinger operators in the diagonal part and small interactions of order $h^\nu$ in the off-diagonal part, 
where $h$ is a semiclassical parameter and $\nu$ is a constant larger than $1/2$.
We study the absence of resonance near a non-trapping energy for both Schr\"odinger operators in the presence of crossings of their potentials. 
The width of resonances is estimated from below by $Mh\log(1/h)$ and the coefficient $M$ is given in terms of the directed cycles of the generalized bicharacteristics induced by two Hamiltonians.
\end{abstract}

\section{Introduction}

We are interested in the resonance free domain for the semiclassical $2\times 2$ matrix  Schr\"odinger operator
\begin{equation}
\begin{aligned}
P(h) &=
\begin{pmatrix}
 P_1 & h^\nu W\\
h^\nu W^* & P_2
\end{pmatrix},
\end{aligned}
\end{equation}
where $$P_j=-h^2\frac{d^2}{dx^2} +V_j(x), \,x\in \mathbb R,\ ( j=1,2). 
$$
Here  
$h>0$ denotes the semiclassical parameter,  $\nu>1/2$, 
$W=W(x,hD_x)$ is a first-order semiclassical differential operator
and $W^*$ denotes  its formal adjoint.
Such  operator appears in the Born-Oppenheimer approximation of
molecules, after reduction to an effective Hamiltonian (see e.g. \cite{KMSW}).

For each semiclassical Schr\"odinger operator  $P_j$ $(j=1,2)$  with $C^\infty$ potential $V_j(x)$,  it is well known that there are no resonances with imaginary part of order $h\log(1/h)$ 
around an energy level $E_0$ satisfying the non-trapping condition  (see \cite{Ma2,SjZw}). 
We recall that an energy $E_0$ is said to be non-trapping if for all compact $K\subset p_j^{-1}(E_0)$ there exists $T_K>0$ such that
\begin{equation}
\label{NT} (x,\xi)\in K\Longrightarrow  {\rm exp}(tH_{p_j})(x,\xi) \not \in K,\ |t|>T_K.
\end{equation}
Here, the function $p_j(x,\xi)=\xi^2+V_j(x)$ defined on the phase space $\R_x\times \R_\xi$ is the symbol of $P_j$ in the sense of Eq. \eqref{WeylQ},  $ H_{p_j}=2\xi \p_x - V'_j(x) \p_\xi$
is the  Hamiltonian vector field, and  $ {\rm exp}(tH_{p_j})(x,\xi)$ is  the corresponding Hamiltonian flow.
It is well-known that (see \cite{GeMa}) the non-trapping condition \eqref{NT} is equivalent to the existence of an escape function $G_j(x,\xi)$  in a neighborhood  of $p_j^{-1}(E_0)$, that is,
\begin{align}\label{NTC}
\exists G_j \in C^\infty(\mathbb R^2;\mathbb   R); \, \, H_{p_j}(G_j)\geq \dl>0 \text { for }  \vert p_j(x,\xi)-E_0\vert \leq \epsilon.
\end{align}

In this article, we consider a non-trapping energy for both $P_1$ and $P_2$.
Under certain conditions, the operator $P$ has the same resonance free domain as in the scalar case. First, under the gap condition (i.e.,  ${\inf_{x\in \mathbb R}}(V_2(x)-V_1(x))>0$), $P(h)$ is unitarily equivalent to a diagonal operator (see  \cite{Sj2}), and we can apply the Martinez's approach \cite{Ma2} to the system. Next, if there exists a common escape function $G$ such that \eqref{NTC} holds for both $p_1$ and $p_2$, then the Sj\"ostrand-Zworski's approach \cite{SjZw} can also be applied to the system $P(h)$ (see \cite{As,ADF}). 
Notice that the gap condition implies the existence of a common escape function $G$. 

\begin{figure}
\centering
\includegraphics[bb=0 0 860 275, width=12cm]{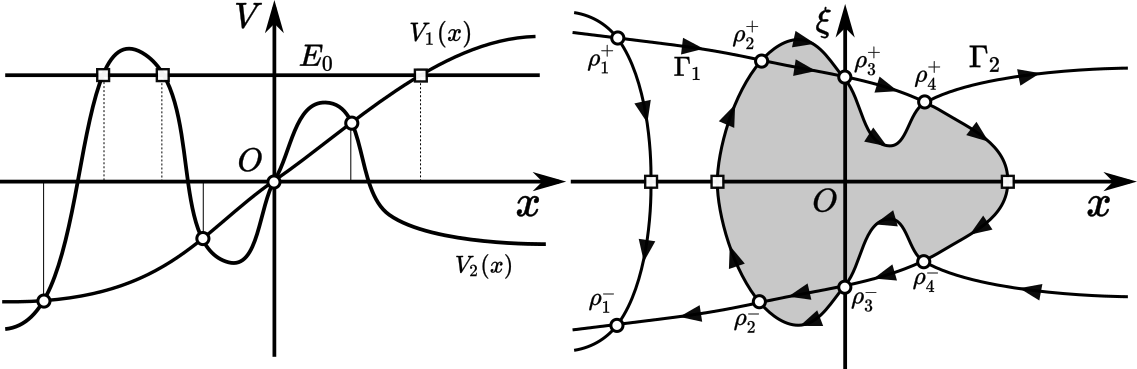}
\caption{$E_0,\,V_1,\,V_2$ satisfying Assumptions (A1)--(A4) and $\vec{\mc{C}}(E_0)\neq\emptyset$.}
\label{Fig:target}
\end{figure}

However, the situation becomes more complicated when the characteristic sets $p_1^{-1}(E_0)$ and $p_2^{-1}(E_0)$ cross each other, that is, when the set of crossing points $\Gamma_c=\Gamma_c(E_0):=p_1^{-1}(E_0)\cap p_2^{-1}(E_0)$ is not empty. 
In particular,
if they form a ``directed cycle" (see Figures 
\ref{Fig:target} and page 5), there is no common escape function.
It is natural to expect to find resonances of $P(h)$ with imaginary part of order $h\log(1/h)$. 

\begin{figure}
\centering
\includegraphics[bb=0 0 880 270, width=11cm]{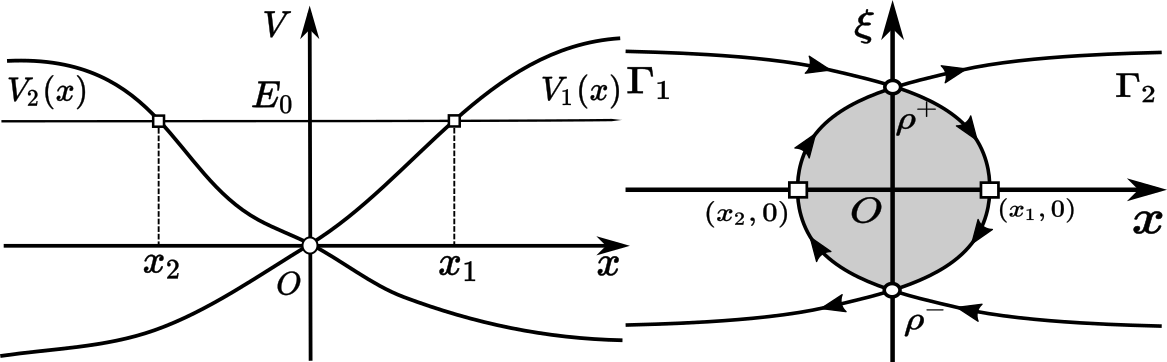}
\caption{A case with a directed cycle (Case T)}
\label{Fig:Lemma1}
\includegraphics[bb=0 0 880 270, width=11cm]{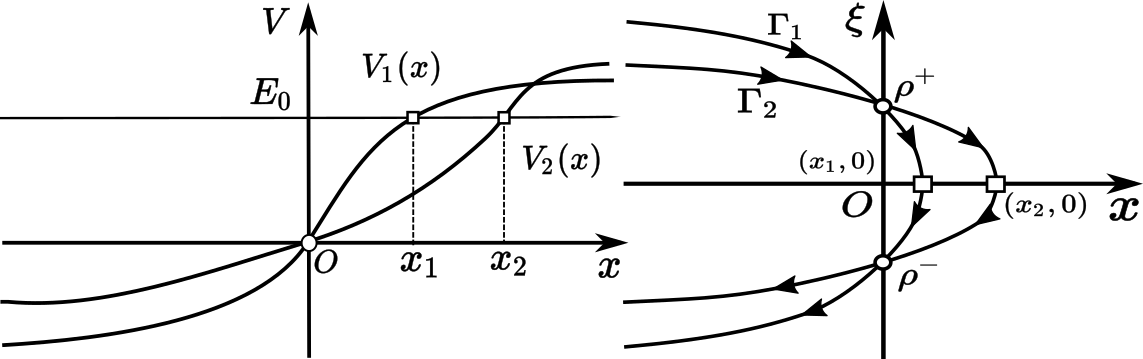}
\caption{A case without any directed cycles (Case N)}
\label{Fig:Lemma2}
\end{figure}

In this paper, we study the absence of resonances for the system $P(h)$ in the presence of finitely many crossing points below $E_0$ especially when there is no common escape function (see Figure \ref{Fig:target}). 
In Theorem \ref{resonancefreedomain}, we obtain the resonance free domain of size $Mh\log(1/h)$ with $M:=(2\nu-1)T^{-1}$. 
Here, $T=T(E_0)$ stands for the time which a classical particle spends traveling along a directed cycle (see page 5).

We emphasize that the proof of Theorem \ref{resonancefreedomain} relies on microlocal argument independent of both Martinez \cite{Ma2} and Sj\"ostrand-Zworski \cite{SjZw}. 
Our argument can be briefly summarized as follows : 
We construct microlocal WKB solutions supported on the characteristic set $\Gamma=\Gamma(E_0):=p_1^{-1}(E_0)\cup p_2^{-1}(E_0)$ except at crossing points and turning points, 
and we connect them by connection formulae. 
We prove that if $|\im E|<Mh\log(1/h)$ with $M$ mentioned above, then the microlocal WKB solution becomes smaller after a connection around certain directed cycle, which is a contradiction.

These microlocal connection formulae were established in \cite{AsFu,FMW3} using normal form in the spirit of \cite{CdvPa,HeSj2,Sj1} (see also \cite{BFRZ1, BFRZ2, Cdv}). 

Finally, we note that, in the case $\nu>1$, we can prove a resonance free domain of order $M'h\log(1/h)$ for some constant $M'>0$ only by a resolvent estimate (see Remark \ref{ResolventEstimate} and \ref{RE}). 
This method works even in the multidimensional case, but does not give the precise value of $M'$.

This work is motivated by the serie of works \cite{Ash,Ba,FMW1,FMW3,Ma1,Na}
where the width of resonances of the same  (or more general) operator as $P$  are studied
for different potentials and different energy levels (see also \cite{AsFu,Pe} for studies concerning eigenvalue splitting for systems).
In particular, in \cite{FMW3}, the asymptotics of
the resonances when the corresponding classical dynamical system has a directed cycle was given. 
The quantum effect of the crossing appears in the asymptotics of the resonances, but the principal part is governed by the Bohr-Sommerfeld quantization rule associated with
the simple well of one of the potentials. In our work,
there is no cycle made only by a single potential, and the effect of crossing appears in the constant $M$.
Finally, it is worth mentioning that in the papers stated above, they only consider weaker interaction than in our case, i.e., the case $\nu=1$ (or $\nu=2$ in \cite{Pe}).

\section{Main Result}\label{AssRes}
We study the absence of resonances in the semiclassical limit $h\to 0_+$ in a neighborhood of an energy $E_0\in\R$.
Throughout this paper, we assume the following (A1)--(A4) :

\vspace{0.2cm}
\noindent
\textbf{(A1)}
For $j=1,2$, $V_j(x)$ is a real-valued smooth function on $\R$, satisfying following conditions:
\begin{enumerate}
\item It extends to a holomorphic function
in an anguler complex domain near infinity
$$
{\mathcal S}=
\bigl\{x\in\C\, ; \,  |\im\, x|<(\tan\theta_0)|\re \,x|,\ |\re\, x|> R_0
\bigr\}
$$
for some constants $0<\theta_0<\pi/2,$ $R_0>0$.
\item It admits limits as $\re x\to\pm\infty$ in $\cS$, and
$$
v_j^\pm:=\lim_{{\re\,x\to \pm\infty}\atop{x\in {\mathcal S}}}V_j(x)\neq E_0.
$$
\item The energy $E_0$ is non-trapping for $P_j(h)$, i.e., the condition \eqref{NT} (or equivalently \eqref{NTC}) holds. 
\end{enumerate}

\noindent
\textbf{(A2)}
The operator $W(x,hD_x)$ has the form
$$
W(x,hD_x)=r_0(x)+ir_1(x)hD_x,\quad D_x:=-i\frac{d}{dx}
$$
where $r_0$ and $r_1$ are bounded smooth real-valued functions on $\R$
and extend to holomorphic functions in ${\mathcal S}$.


\vskip 0.3cm
In this situation, $P$ is essentially self-adjoint with domain $H^2(\R)\oplus H^2(\R)$, and the resonances of $P$ can be defined, e.g., 
as the values $E\in\C_-=\{\im z<0\}$ such that the equation $Pu=Eu$ has a non trivial outgoing solution $u$, that is, a non identically vanishing solution such that, for some $\theta \in]0,\theta_0]$, the function $u\circ\zeta_\theta$ is in $L^2(\R)\oplus L^2(\R)$ 
where $\zeta_\theta\in C^\infty(\R)$ satisfies $\zeta_\theta(x)= x$ for $|x|\le R_0$ and $\zeta_\theta(x)=xe^{i\theta}$ for $|x|\ge2R_0$ (see, e.g., \cite{AgCo,DyZw, ReSi}). 
Equivalently, the resonances are the eigenvalues of the operator $P_\theta=U_\theta P U_\theta^{-1}$ acting on $L^2(\R_\theta)\oplus L^2(\R_\theta)$, 
where $\R_\theta=\zeta_\theta(\R)$ is a complex dilation of $\R$ and $U_\theta u=|\zeta_\theta'(x)|^{1/2}u\circ\zeta_\theta$ (see, e.g., \cite{HeMa}). 
Note that there is no essential spectrum in some $h$-independent complex neighborhood of $E_0$ due to the assumption $v_j^\pm\neq E_0$. 
We denote by ${\rm Res}(P)$ the set of these resonances.

\vspace{0.2cm}
We denote the characteristic set of $P_j$ by 
$$
\Gamma_j=\Gamma_j(E_0):=p_j^{-1}(E_0)=\{(x,\xi)\in\R^2\, ; \, \xi^2 +V_j(x)=E_0\} \quad(j=1,2).
$$
We define the set of crossing points in the phase space by $\Gamma_c=\Gamma_c(E_0):=\Gamma_1\cap\Gamma_2$. 
Let $w(x,\xi)=r_0(x)+ir_1(x)\xi$ be the principal symbol of the operator $W$.

\vspace{0.2cm}
\noindent
\textbf{(A3)}
The set $\Gamma_c(E_0)$ is finite. $\Gamma_1$ and $\Gamma_2$ intersect transversally at every point of $\Gamma_c$, and they are not asymptotic to each other (i.e., $v_1^+\neq v_2^+$ and $v_1^-\neq v_2^-$ holds). 

\vspace{0.2cm}
\noindent
\textbf{(A4)}
The operator $W$ is microlocally elliptic at each point of $\Gamma_c$ (i.e., $w(\rho)\neq0$ holds for all $\rho\in\Gamma_c$) .

\vspace{0.2cm}
We introduce a graph structure on $\Gamma$.  
Let $\Gamma_t=\Gamma_t(E_0):=\{(x,0);\,x\in(V_1^{-1}(E_0)\cup V_2^{-1}(E_0))\}$ be the set of turning points, and 
$$
\mathcal{V}=\mc{V}(E_0):=\Gamma_c\cup \Gamma_t,\quad\mc{V}_j:=\mc{V}\cap\Gamma_j\quad\mbox{(vertices)}.
$$
Then, the set $\Gamma\setminus\mc{V}$ can be uniquely written as a disjoint union of connected components $\Gamma\setminus\mc{V}=\bigsqcup\widetilde{\mc{E}}$, 
and we put
\begin{align*}
&\mathcal{E}=\mc{E}(E_0):=\{\gamma\in\widetilde{\mc{E}};\, \gamma:\mbox{bounded}\},
\quad\mathcal{E}_j=\mc{E}_j(E_0):=\{\gamma\in\mc{E};\, \gamma\subset\Gamma_j\}\quad\mbox{(edges)}.
\end{align*}
Then we have $\mc{E}=\mc{E}_1\sqcup\mc{E}_2$, and we can see the pairs $\mc{G}:=(\mathcal{E},\mathcal{V})$ and $\mc{G}_j:=(\mc{E}_j,\mc{V}_j)$ as graphs by identifying each element of $\mc{E}$ with the pair of its endpoints. 
We see that some cycle may appear in the graph $\mc{G}$ when $\Gamma_c\neq\emptyset$.
Next, we give an orientation to the graph by Hamiltonian vector field. 
In other words, let the time function $t$ given by the following path integral
\begin{align}\label{Timefn1}
t(\gamma)
:=\int_{\gamma }\frac{dx}{2\xi} \quad(\gamma\in\mc{E})
\end{align}
be positive for all $\gamma\in\mc{E}$.
Then, one has 
$$
\exp(t(\gamma)H_{p_j})(\gamma_-)=\gamma_+\ \mbox{ for any }\gamma\in\mc{E}_j,
$$ 
where $\gamma_-$ and $\gamma_+$ stand for the initial and the terminal points of $\gamma$ respectively after orientation ($\gamma_\pm\in\mc{V}$). 
This means that $t(\gamma)$ is the time which a classical particle spends traveling along $\gamma$ from $\gamma_- $ to $\gamma_+$. 
We can extend the time function $t$ to the paths (consecutive sequence of edges) on $\mc{G}$ by
\begin{align}\label{Timefn2}
t(\gamma):=\sum_{k=1}^m t(\gamma_k),\quad
\left(\gamma=\bigcup_{k=1}^m \gamma_k: \mbox{a path},\ \gamma_k\in\mc{E}\right).
\end{align}
We denote $\vec{\mc{C}}(E_0)$ the set of directed cycles, 
and $\vec{\mc{C}}_{\ope{pr}}(E_0)$ the set of directed cycles consisting of two paths belonging to $\Gamma_1$ and $\Gamma_2$, respectively
$$
\vec{\mc{C}}_{\ope{pr}}(E_0):=\{\gamma=\gamma_1\cup\gamma_2\in\vec{\mc{C}}(E_0)\,; \, \gamma_j\mbox{ is a path in }\mc{G}_j\ (j=1,2)\}.
$$
Notice that $\vec{\mc{C}}_{\ope{pr}}(E_0)$ is finite, and not empty whenever $\vec{\mc{C}}(E_0)\neq\emptyset$.  

\begin{theorem}\label{resonancefreedomain}
Assume (A1)--(A4). 
Then, there exists a positive constant $M$ such that
\begin{align}
\label{2.3}
\ope{Res}(P(h))\cap \{z\in\C_-;\, |z-E_0|<Mh\log(1/h)\}=\emptyset
\end{align}
holds for $h$ small enough. 
Here, the constant $M>0$ is given by \begin{align}\label{boundaryT}
M=\frac{2\nu-1}{T(E_0)}\quad\mbox{with}\quad 
T(E_0)=\max_{\gamma\in\vec{\mc{C}}_{\ope{pr}}(E_0)}t(\gamma),
\end{align}
if $\vec{\mc{C}}(E_0)\neq\emptyset$.
Otherwise, it is arbitrary.
\end{theorem}

\begin{remark}
\begin{enumerate}
\item If $\vec{\mc{C}}(E_0)=\emptyset$, 
we have a common escape function for both $P_1$ and $P_2$, 
and we can apply Sj\"ostrand-Zworski's approach \cite{SjZw} to show \eqref{2.3} with arbitrary $M$ (see \cite{As,ADF}). 
Our method gives a different proof.
\item Since the set of energies $E$ satisfying Assumptions (A1)--(A4) is open, 
there exists some interval $I$ such that our result holds uniformly in $\re E\in I$.
\end{enumerate}
\end{remark}
\begin{remark}\label{ResolventEstimate}
Under Assumptions (A1) and (A2), there exist positive $(h,\nu)$-independent constants $M_1$ and $h_0$ such that there are no resonances of $P(h)$ in $\{z\in\C_-;\, |z-E_0|<Mh\log(1/h)\}$ for all $h\in]0,h_0]$ with $M=M_1(\nu-1)$ (only for $\nu>1$). 
This holds for $d$-dimensional case, that is, 
\begin{align}\label{multiD}
P_j=-h^2\bigtriangleup+V_j(x)\quad\mbox{and}\quad
W=r_0(x)+\sum_{k=1}^dr_{1,k}(x)h\frac{d}{dx_k},\quad(x\in\R^d)
\end{align}
with $r_0,\,r_{1,k},\, V_j\in C^\infty(\R^d;\R)$ 
extend analytically to $\cS^d$ ($k=1,\ldots,d$, $j=1,2$), $V_j$ admits limit $v_j^\infty\neq E_0$ as $|x|\to+\infty$ in $\cS^d$. 
We give the proof in \ref{RE}.
\end{remark}

\section{Proof of the Theorem}
In this section, we give a proof of Theorem \ref{resonancefreedomain}.
\subsection{Semiclassical and microlocal terminologies}
Here, we recall some basic notions of semiclassical and mirolocal analysis, referring to the books \cite{DiSj, Ma3, Zw} for more details.
Let $S^m$ 
be the space of symbols
$$
S^m:=\left\{a\in C^\infty(\R^2;\C)\,; \, 
\left|\p_x^k \p_\xi^l a\right|\le C_{k,l}(1+\xi^2)^{(m-l)/2}
\mbox{ for all }k,l\in\N\right\}
$$
with $\N=\{0,1,2,\ldots\}$.
The $h$-pseudodifferential operator corresponding to a symbol $a\in S^m$ denoted $\wey{a}(x,hD)$ is defined on Sobolev space $H^m(\R)$ by
\begin{align}\label{WeylQ}
\wey{a}(x,hD)u(x)
:=\frac{1}{2\pi h}\int_{\R^2}e^{i(x-y)\xi/h}a\left(\frac{x+y}{2},\xi\right)u(y)dyd\xi.
\end{align}
Let $U\subset \R^2$ and $f=f(x;h)\in L^2(\R)$ with $\|f\|_{L^2}\le 1$. 
We say $f$ is microlocally $0$ in $U$ and write
$$
f(x;h)\sim0\mbox{ microlocally in }U
$$
if there exists a symbol $\chi(x,\xi)\in S^0(\R^2)$ with $\chi(x,\xi)\neq0$ for all $(x,\xi)\in U$  
such that 
$$
\wey{\chi}(x,hD)f(x;h)=\ord(h^\infty)\mbox{ uniformly in }\R^2.
$$
We define $\ope{WF}_h(f)$ the semiclassical wave front set of $f$ as the set of all points of $\R^2$ such that $f$ is not microlocally $0$. 
The followings are well-known (see e.g. \cite[Theorem 12.5]{Zw}) : 
If $f$ satisfies $a^w(x,hD)f\sim0$ microlocally in $U$, we have $\ope{WF}_h(f)\subset \{a=0\}$. 
Moreover, if $\p a\neq0$ on $\{a=0\}$, the semiclassical wave front set $\ope{WF}_h(f)$ is invariant under the Hamiltonian flow of $a$. 
According to them, in our problem, the semiclassical wave front set of a resonant state $u$ is a subset of the characteristic set $\Gamma$, and the space of microlocal solutions on each edge $\gamma\in\mc{E}$ is one-dimensional. 

Let $\Lambda$ be the set defined by
$
\Lambda=\{k+l\nu\, ; \, k,l\in\N\}.
$
We say that an $h$-dependent constant $\mu(h)$ has an asymptotic expansion 
\begin{align}\label{Asymptotic1}
\mu(h)\sim\sum_{\lambda\in\Lambda} h^\lambda\mu_\lambda
\end{align}
if for every $\lambda_0\in\Lambda$, there exists a constant $C_{\lambda_0}$ 
such that 
\begin{align}\label{Asymptotic2}
\biggl|\mu(h)-\sum_{\lambda\in\Lambda,\lambda< \lambda_0}h^\lambda\mu_\lambda\biggr|<C_{\lambda_0}h^{\lambda_0}
\end{align}
holds for small $h>0$. Here, $\mu,\,\mu_\lambda$ can depend on some other parameters.

\subsection{One point crossing case}\label{proof} 
For simplicity, in this subsection, we give the proof of our Theorem in the case of one crossing point (see Figures \ref{Fig:Lemma1},\ref{Fig:Lemma2}), and we send the reader to subsection \ref{finitelymany} for the proof in the general case. 
We take $M'>0$ arbitrary, and assume that throughout this section $E\in\cB_h(M')=\{z\in\C_-;\,|z-E_0|<M'h\log(1/h)\}$. 
We assume that
\begin{align*}
\{V_1=V_2\le E_0\}=\{0\},\quad
V_1(0)=V_2(0)=0,\quad 
V_1'(0)>V_2'(0),\quad
v_1^-<E_0<v_1^+,
\end{align*}
and 
there exist $x_1,\,x_2\in\R$ such that $V_j^{-1}(E_0)=\{x_j\}$ holds for both $j=1,2$. 
In this situation, existence and non-existence of directed cycle is equivalent to
$$
v_2^->E_0>v_2^+\ \mbox{(Case T)},\quad 
\mbox{and}\quad v_2^-<E_0<v_2^+\ \mbox{(Case N)},
$$ 
respectively (also equivalent to $x_2<0<x_1$ and $0<x_1<x_2$), 
and the sets $\Gamma_1$ and $\Gamma_2$ cross at two points $\rho^\pm (E_0):=(0,\pm\sqrt{E_0})$. 
For  $j=1,2$, we also put (see Figures \ref{Fig:Lemma1} and \ref{Fig:Lemma2}),
\begin{align*}
&\Gamma_{j,L}=\{(x,\xi)\in \Gamma_j \, ; \, x<0\},\quad
\gamma_{j,L}^\pm=\{(x,\xi)\in \Gamma_j\, ; \, x<0,\pm\xi>0\},\\
&\Gamma_{j,R}=\{(x,\xi)\in \Gamma_j\, ; \, x>0\},\quad
\gamma_{j,R}^\pm=\{(x,\xi)\in \Gamma_j\, ; \, x>0,\pm\xi>0\}.
\end{align*}
By these notation, we have $\Gamma_c=\{\rho^\pm\},$ $\Gamma_t=\{x_1,x_2\}$, $\til{\mc{E}}=\{\gamma_{j,S}^\pm;\,(j,S)\in\{1,2\}\times\{L,R\}\}$.

\subsubsection{Review of microlocal connection formulae}
In this section, we recall the microlocal connection formulae at crossing and turning points from \cite[Sec. 5]{FMW3} and \cite[Sec. 4]{AsFu}. Notice that, in these works, the considered energy $E$ is complex with imaginary part of order $h$ in \cite{FMW3} and real in \cite{AsFu}. However, in our case, we deal with complex energies with imaginary part of order $h\log(1/h)$. Thus, some terms which were negligible as $h\to0_+$ have to be controlled in our case.
Also, \cite{FMW3} suppose that $V_1$ and $V_2$ are analytic, and we apply normal form of \cite{Sj1} same as in \cite{AsFu}. 
Note that the definition of asymptotic series is extended by Eqs. \eqref{Asymptotic1}, \eqref{Asymptotic2}, and we can treat the parameter $\nu$.

For each $(j,S)\in\{1,2\}\times\{L,R\}$, there exist functions $f_{j,S}^\pm=f_{j,S}^\pm(x;h)$ such that
\begin{align}
(P-E)f_{j,S}^\pm\sim 0,
\quad f_{j,S}^\pm\sim 
e^{\pm i\phi_j^0(x)/h}
(a_j^{\pm}, b_j^{\pm})^T
\,\,\,{\rm microlocally\, on}\,\,
\gamma_{j,S}^\pm.
\end{align}
Here, the phase function $\phi_j^0(x)$ is given by
$$
\phi_j^0(x):= \int_0^x \left(\sqrt{ E_0-V_j(t)}+\frac{E-E_0}{2\sqrt{E_0-V_j(t)}}\right) \,dt ,
$$
and 
$a_j^\pm= a_j^\pm(x;h)\sim \sum_{\lambda\in\Lambda}h^{\lambda} a_{j,\lambda}^\pm(x)$ and 
$b_j^\pm=b_j^\pm (x;h)\sim \sum_{\lambda\in\Lambda}h^{\lambda} b_{j,\lambda}^\pm(x)$ are 
symbols whose first coefficients are given by 
\begin{align}
\label{coeff1}
\left\{
\begin{aligned}
&a_{1,0}^\pm(x)
=\frac{1}{(E_0-V_1(x))^{1/4}} \, ,\\
&b_{1,\nu}^\pm(x)
=\frac{r_0(x)\mp ir_1(x)\sqrt{E_0-V_1(x)}}{(V_1(x)-V_2(x))(E_0-V_1(x))^{1/4}} .
\end{aligned}
\right.
\end{align}
Similar computation on $\gamma_{2,S}^\pm$ ($S=L,R$) leads $a_{2,\lambda}^\pm\equiv0$ ($\lambda<\nu$), and
\begin{align}
\label{coeff2}
\left\{
\begin{aligned}
&b_{2,0}^\pm(x)
=\frac{1}{(E_0-V_2(x))^{1/4}} \, ,\\
&a_{2,\nu}^\pm(x)
=\frac{-r_0(x)\mp ir_1(x)\sqrt{E_0-V_2(x)}}{(V_1(x)-V_2(x))(E_0-V_2(x))^{1/4}} .
\end{aligned}
\right.
\end{align}
In particular, $b_{1,\lambda}^\pm\equiv a_{2,\lambda}^\pm\equiv0$ holds for all $\lambda<\nu$.

Since the space of microlocal solutions on each edge is one-dimmensional, 
if a function $u$ is a resonant state corresponding to $E$, 
there exist constants $\alpha_{j,S}^\pm=\alpha_{j,S}^\pm(h)$ for $(j,S)\in\{1,2\}\times\{L,R\}$ such that
$$
u\sim  \alpha_{j,S}^\pm f_{j,S}^\pm
\,\,\,{\rm microlocally\, on}\,\,
\gamma_{j,S}^\pm.
$$
Then it holds that
\begin{align}\label{TMatpm}
&
\begin{pmatrix}
\alpha^+_{1,R}\\
\alpha^+_{2,R}\\
\end{pmatrix}
=
\begin{pmatrix}
\tau_{1,1}^+	&\tau_{1,2}^+\\
\tau_{2,1}^+	&\tau_{2,2}^+\\
\end{pmatrix}
\begin{pmatrix}
\alpha^+_{1,L}\\
\alpha^+_{2,L}\\
\end{pmatrix},\quad
\begin{pmatrix}
\alpha^-_{2,L}\\
\alpha^-_{1,L}\\
\end{pmatrix}
=
\begin{pmatrix}
\tau_{1,1}^-	&\tau_{1,2}^-\\
\tau_{2,1}^-	&\tau_{2,2}^-\\
\end{pmatrix}
\begin{pmatrix}
\alpha^-_{2,R}\\
\alpha^-_{1,R}\\
\end{pmatrix},
\end{align}
and 
\begin{align}\label{MaslovTN}
\left\{
\begin{aligned}
&\alpha_{1,R}^-=m_1\alpha_{1,R}^+,\mbox{ and }\alpha_{2,L}^-=m_2\alpha_{2,L}^+	&\mbox{(Case T)}\\
&\alpha_{j,R}^-=m_j\alpha_{j,R}^- 	&\mbox{(Case N)}.
\end{aligned}\right.
\end{align}
where the constants $\tau_{j,k}^\pm=\tau_{j,k}^\pm(E;h)$ and $m_j=m_j(E;h)$ behave
\begin{align}
&
\begin{aligned}
&\tau_{1,1}^\pm=1+\ord(h^{\nu'}),			&\tau_{1,2}^\pm&=-ih^{\nu-\omega_1^\pm-\frac{1}{2}} (\sigma^\pm +\ord(h^{\nu'})),\\
&\tau_{2,1}^\pm=-ih^{\nu-\omega_2^\pm-\frac{1}{2}}(\overline{\sigma}^\pm +\ord(h^{\nu'})) ,	&\tau_{2,2}^\pm&=1+\ord(h^{\nu'}),
\end{aligned}\\
&m_j=
(-i+\ord(h^{\nu'}))e^{2i(S_j+(E-E_0)T_j)/h},
\end{align}
as $h\to0_+$ uniformly for $E\in\cB_h(M')$.
Here, the constants $\nu',\,\omega_j^\pm,\,\sigma^\pm,\,S_j,\,T_j$ are given by
\begin{align}\label{constants}
&\nu'=\min\{1,\nu\},\quad
\omega_j^\pm:=i\mu_j^\pm h^{-1+2\nu},\quad
\sigma^\pm:= e^{\pm i\pi /4}E_0^{-1/4}w(\rho^\pm)\sqrt{\frac{\pi}{v_1-v_2}},\nonumber\\
&S_j=\int_0^{x_j}\sqrt{E_0-V_j(x)}\,dx,\quad T_j=\int_0^{x_j}\frac{dx}{\sqrt{E_0-V_j(x)}},
\end{align}
with the constants $\mu_j^\pm=\ord(1)$. Moreover, when $\nu\in\,]1/2,1[$, we have
$$
\mu_j^\pm=\frac{|w(\rho^\pm)|}{2\sqrt{E_0}(v_1-v_2)}+\ord(h^{2(1-\nu)}),
\quad v_j:=V_j'(0).
$$
Notice that 
$\overline{\sigma}^\pm=\sigma^\mp$.

\begin{remark}\label{TMrem}
\begin{enumerate}
\item Eq. \eqref{TMatpm} suggests that 
a ``small difference" between two Schr\"odinger operators ($|V_1'(0)-V_2'(0)|=|\tau_1-\tau_2|$ small),
a ``small energy" ($E_0$ small) 
and a ``big interaction" at $\rho^\pm$ ($|w(\rho^\pm)|$ big and $\nu$ small) 
make it easier to transfer between two charachteristic sets $\Gamma_1$ and $\Gamma_2$.
\item For each $j=1,2$, the value $2S_j$ is the area of region enclosed by $\xi$-axis and $\Gamma_j$ in the phase space, 
and the value $2T_j$ satisfies $T_j=t(\gamma_{j,S}^+\cup\gamma_{j,S}^-)$ by chosing suitable $S=L,R$.
\end{enumerate}
\end{remark}

\subsubsection{Contradiction argument}\label{Contradiction}
We give a proof for the case with only one crossing point by using the previous connection formulae. 
Suppose that $E\in \cB_h(M')$ is a resonance of $P$, 
and that $u(x;h)=(u_1(x;h),u_2(x,h))^T$ is a correspoding resonant state. 
Then, $u$ is microlocally $0$ in incoming region, that is, 
\begin{align}\label{InCom}
\left\{
\begin{aligned}
&u\sim0 \mbox{ microlocally on }\gamma_{1,L}^+\cup\gamma_{2,R}^-&\mbox{(Case T)}\\
&u\sim0 \mbox{ microlocally on }\gamma_{1,L}^+\cup\gamma_{2,L}^+&\mbox{(Case N)}.
\end{aligned}
\right. 
\end{align}
For Case T, 
we can assume that $\alpha_{2,L}^+=1$, i.e., 
$
u\sim f_{2,L}^+\mbox{ microlocally on }\gamma_{2,L}^+.
$
Then by \eqref{TMatpm} and by $u\sim0$ microlocally on $\gamma_{1,L}^+$, 
we have 
$$
u\sim \tau_{2,2}^+ f_{1,R}^+\mbox{ microlocally on }\gamma_{1,R}^+.
$$
Similarly, by Eqs. \eqref{TMatpm}, \eqref{MaslovTN} and by $u\sim0$ microlocally on $\gamma_{2,R}^-$, 
we have 
$$
u\sim C(E;h)f_{2,L}^+\mbox{ microlocally on }\gamma_{2,L}^+,
$$
with $C(E;h):=m_2\tau_{1,2}^-m_1 \tau_{1,2}^+$. 
Here, we have
$$
\begin{aligned}
C(E;h)
=
\left(|\sigma^+|^2 +\ord(h^{\nu'})\right)
h^{2\nu-1-(\omega_1^+ +\omega_2^-)}
\exp\left(\frac{i}{h}(S+(E-E_0)T)\right),
\end{aligned}
$$
where $S=S(E_0):=2(S_1+S_2)$ is the area of the domain bounded by the directed cycle 
$\gamma:=(\gamma_{1,R}^+\cup\gamma_{1,R}^-)\cup(\gamma_{2,L}^-\cup\gamma_{2,L}^+)\in\vec{\mc{C}}_{\ope{pr}}(E_0)$ in the phase space,
and 
$$
T=T(E_0)=2(T_1+T_2)=t(\gamma)
$$ is the time required for a classical particle to travel along $\gamma$. 
(Note that $|\sigma^+|^2=\sigma^+\sigma^-$.) 
In particular, 
\begin{align}
\label{order of all}
C(E;h)
=\ord(h^{2\nu-1}e^{T(E_0)|\im E|/h})
\end{align}
holds. We see immediately that if $E\in\cB_h(M)$, $|C(E;h)|\ll1$ for small $h$. Because $E\in\cB_h(M)$ implies that $|\im E|<\frac{2\nu-1}{T(E_0)}h\log(1/h)$.
This leads to a contradiction. 
Hence, we deduce that there are no resonances of $P(h)$ in $\cB_h(M)$.

For Case N, we have $u\sim0$ microlocally on $\gamma_{1,L}^+\cup \gamma_{2,L}^+$ by \eqref{InCom}.
So by \eqref{TMatpm}, we also have 
$$
u\sim 0\mbox{ microlocally on }\gamma_{1,R}^+\cup\gamma_{2,R}^+.
$$
Similarly, Eqs. \eqref{MaslovTN}, \eqref{TMatpm} imply that $u\sim 0$ microlocally $\gamma_{1,R}^-\cup\gamma_{2,R}^-$ and $\gamma_{1,L}^-\cup\gamma_{2,L}^-$, respectively.
This leads that $u\equiv \ord(h^\infty)$ and contradicts with the fact that $u$ is a resonant state (non identically vanishing).

\subsection{The case with finitely many crossing points}\label{finitelymany}

Here, we consider the general case. Arguments are essentially same as the one point crossing case. 
If $\vec{\mc{C}}(E_0)\neq\emptyset$, using connection formulae, we obtain a necessary condition $C(E;h)=1$ by an almost same way as Case T above with 
\begin{align}\label{GeneC}
C(E;h)=(1+\ord(h^{\nu'}))\sum_{\gamma\in\vec{\mc{C}_2}(E_0)} C_\gamma(E;h).
\end{align}
Notice that for each directed cycle $\gamma\in\vec{\mc{C}_2}(E_0)$, there are two paths $\gamma_1$ and $\gamma_2$ of $\mc{G}_1$ and $\mc{G}_2$ such that $\gamma=\gamma_1\cup\gamma_2$, and we denote the common initial point of $\gamma_j$ (i.e., terminal point of the other) by $\rho_j$. 
Then the constants $C_\gamma(E;h)$ in \eqref{GeneC} are given by
\begin{align*}
\left(\sigma_{\rho_1}{\sigma}_{\rho_2}+\ord(h^{\nu'})\right)
h^{2\nu-1-(\omega_{\rho_1}+\omega_{\rho_2})}
\exp\left(\frac{i}{h}(S_\gamma+(E-E_0)T_\gamma)\right),
\end{align*}
where the constats $\sigma_{\rho_j},\,\omega_{\rho_j}$ are given similarly as \eqref{constants}. 
In particular,
\begin{align*}
S_\gamma=\int_\gamma \xi dx,\quad T_\gamma=\int_\gamma\frac{dx}{2\xi}=t(\gamma).
\end{align*}
If $E\in\cB_h(M)$, $C_\gamma(E;h)\to0$ as $h\to0_+$ for all $\gamma\in\vec{\mc{C}}_{\ope{pr}}(E_0)$. 
Consequently, we see that $C(E;h)\to0$ as $h\to0_+$. This leads to a contradiction.

When $\vec{\mc{C}}(E_0)=\emptyset$, the argument is exactly the same as Case N above, that is, we obtain $u\equiv\ord(h^\infty)$ after connecting.

\section*{Acknowledgements}
Special thanks should be given to Maher Zerzeri for many stimulating discussions during the author's visit to University Paris 13 and to Ritsumeikan University for the financial support. 
The author would like to thank Mouez Dimassi for the important remarks in this research. 
The author is grateful to the organizers of ``Summer Northwestern Analysis Program" held at Northwestern University in August, 2019 where he enjoyed lectures and discussions. 

\appendix
\section{Resolvent Estimate}\label{RE}
Here, we give the proof of Remark \ref{ResolventEstimate}. 
The operators $P_j,W$ are given by Eq. \eqref{multiD}, and the definition of escape function extends naturally to $d$-dimensional case. 
Take positive constants $M'>M>0$ (to be determined later) and $E\in\cB_h(M)=\{z\in\C_-;\, |z-E_0|<Mh\log(1/h)\}$. 
By (A1), we have an escape function $G_j(x,\xi)\in C^\infty(\R^{2d};\R)$ for each $P_j$ ($j=1,2$). 
Put 
\begin{align}
P_0:=
\begin{pmatrix}
P_1&0\\
0&P_2\\
\end{pmatrix},
\quad
W_0:=
\begin{pmatrix}
0&W\\
W^*&0\\
\end{pmatrix},\quad
\widetilde{U}:=
\begin{pmatrix}
e^{-\e \til{G}_1^w/h}&0\\
0&e^{-\e \til{G}_2^w/h}
\end{pmatrix},
\end{align}
where 
$\e=M'h\log(1/h)$, $\theta=\tan^{-1}(\e)$, and $\til{G}_j=G_j-\zeta_0(x)\cdot\xi$, 
where $\zeta_0$ is an $\R^d$-valued function $\zeta_0$ satisfying 
$\til{\zeta}_\theta(\R^d)=\zeta_\theta(\R^d)$ with 
$\til{\zeta}_\theta(x):=x+i(\tan\theta)\zeta_0(x).$ 
Then, we have 
\begin{align}\label{dilatedP}
\til{U} (P_\theta-E) \til{U}^{-1}
=\left[I+h^\nu W_{0,\theta}(P_{0,\theta}-E)^{-1}\right](P_{0,\theta}-E),
\end{align} 
where $P_{0,\theta}=\til{U}U_\theta P_0 U_\theta^{-1}\til{U}^{-1},$ $W_{0,\theta}=\til{U}U_\theta W_0 U_\theta^{-1}\til{U}^{-1}$. 
By the approach of Sj\"ostrand-Zworski \cite[page 401]{SjZw}, there exists a constant $c$ such that 
\begin{align}\label{scalarestimate}
\left\|e^{-\e \wey{G_j}/h}U_\theta P_jU_\theta^{-1}e^{\e\wey{G_j}/h}-E_0\right\|\ge  c\e,
\end{align}
and consequently, if $M\le c'M'$ with some $0<c'<c$, we have
\begin{align}
\left\|\left(P_{0,\theta}-E\right)^{-1}\right\|\le \frac{C_1}{h\log(1/h)}.
\end{align}
On the other hand, for the $(1,2)$-element of the operator $W_{0,\theta}$, we have 
\begin{align}
\left\|e^{-\e G_1^w/h}U_\theta WU_\theta^{-1}e^{\e G_2^w/h}\right\|
\le C_2 e^{M'\|G_2-G_1\|_\infty\log (1/h)}
=C_2h^{-M'\|G_2-G_1\|_\infty}.
\end{align}
Here, we can assume that $G_j(x,\xi)=x\cdot\xi$ for $|(x,\xi)|$ large enough, then $G_2-G_1$ is bounded (i.e. $\|G_2-G_1\|_\infty<+\infty$).
Therefore,
\begin{align}
h^\nu\left\|W_{0,\theta}\left(P_{0,\theta}-E\right)^{-1}\right\|
\le C_3\frac{h^{\nu-1-M'\|G_2-G_1\|_\infty}}{\log(1/h)}.
\end{align}
Thus, if $M'\le (\nu-1)\|G_2-G_1\|_\infty^{-1}$, 
the operator $I+h^\nu W_{0,\theta}(P_{0,\theta}-E)^{-1} $ is invertible for small $h$. 
This means that $E$ is an element of the resolvent set of the operator $P_\theta$, hence $E\notin\ope{Res}(P)$.

Kenta Higuchi, Department of Mathematical Sciences, Ritsumeikan University, 1-1-1 Noji-Higashi, Kusatsu, 525-8577,  Japan

\textit{E-mail address}: ra0039vv@ed.ritsumei.ac.jp

\end{document}